\documentclass[12pt]{iopart}

\usepackage{iopams,setstack}
\usepackage{amssymb,amsthm}
\usepackage{graphicx,float}
\usepackage[colorlinks=true,linkcolor=blue,citecolor=blue,urlcolor=blue]{hyperref}
\usepackage{color}
\usepackage{bm} 
\usepackage{verbatim}
\usepackage{fullpage}

\newcommand{\dsp}{\displaystyle}

\makeatletter
\newcommand{\bigintss}{\@ifnextchar_\@bigintsssub\@bigintssnosub}
\def\@bigintsssub_#1{\def\@int@subscript{#1}\@ifnextchar^\@bigintsssubsup\@bigintsssubnosup}
\def\@bigintsssubsup^#1{\mathop{\text{\LARGE$\int_{\text{\normalsize$\scriptstyle\kern-0.25em\@int@subscript$}}^{\text{\normalsize$\scriptstyle#1$}}$}}\nolimits}
\def\@bigintsssubnosup{\mathop{\text{\LARGE$\int_{\text{\normalsize$\scriptstyle\@int@subscript$}}$}}\nolimits}
\def\@bigintssnosub{\@ifnextchar^\@bigintssnosubsup\@bigintssnosubnosup}
\def\@bigintssnosubsup^#1{\mathop{\text{\LARGE$\int^{\text{\normalsize$\scriptstyle#1$}}$}}\nolimits}
\def\@bigintssnosubnosup{\mathop{\text{\LARGE$\int$}}\nolimits}
\makeatother

\newcommand{\be}{\begin{equation}}
\newcommand{\ee}{\end{equation}}

\newcommand{\bes}{\begin{equation*}}
\newcommand{\ees}{\end{equation*}}

\usepackage{tikz}
\usepackage{ulem}

\begin{document}
\title[Iterative algorithms for a non--linear inverse problem in atmospheric lidar]{Iterative algorithms for a non--linear inverse problem in atmospheric lidar}
\author{Giulia Denevi$^{1,3}$, Sara Garbarino$^2$ and Alberto Sorrentino$^{3,4}$}

\address{$^1$ Istituto Italiano di Tecnologia, Genova \\
$^2$ Centre for Medical Image Computing, Department of Computer Science, University College London, London\\
$^3$ Dipartimento di Matematica, Universit\`a di Genova, Genova\\
$^4$ CNR--SPIN, Genova }

\begin{abstract}
We consider the inverse problem of retrieving aerosol extinction coefficients from Raman lidar measurements. 
In this problem the unknown and the data are related through the exponential of a linear operator, the unknown is non--negative and the data follow the Poisson distribution. Standard methods work on the log--transformed data and solve the resulting linear inverse problem, but neglect to take into account the noise statistics.
In this study we show that proper modelling of the noise distribution can improve substantially the quality of the reconstructed extinction profiles.
To achieve this goal, we consider the non--linear inverse problem with non--negativity constraint, and propose two iterative algorithms derived using the Karush-Kuhn-Tucker conditions. 
We validate the algorithms with synthetic and experimental data. As expected, the proposed algorithms out--perform standard methods in terms of sensitivity to noise and reliability of the estimated profile.
\end{abstract}

\section{Introduction}

Lidars are widely used instruments that collect information on remote objects by illuminating them with a laser and then measuring the amount of backscattered light and the corresponding traveling time.
Lidars have many different applications, including ecosystem mapping \cite{lefsky2002lidar}, navigation in urban environments \cite{levinson2007map} and stratigraphic modelling in geology \cite{bellian2005digital}. In this study we are concerned with the atmospheric application of lidar: here the laser is headed towards the sky, and the backscattered light carries information about the optical parameters of the atmosphere, namely extinction and backscattering coefficient, as functions of the height \cite{ansmann1992independent}. These, in turn, can be used to infer the presence and composition of clouds and aerosols at various altitudes, in order to characterize pollution or the atmospheric consequences of disruptive events \cite{karyampudi1999validation,ansmann201016,di2012raman}.

More specifically, in this study we focus on the estimation of the extinction coefficient from so--called \textit{Raman} lidar data \cite{ansmann1990measurement}, where part of the backscattered light arrives at a different wavelength than that emitted by the laser, as a consequence of \textit{inelastic} scattering; this allows us to assume knowledge of the backscattering coefficient, which is unknown for \textit{elastic} scattering. As a consequence, the noiseless lidar equation can be written as

\begin{equation}
P(z) = d(z) e^{\displaystyle - \int_0^z \alpha(x)dx} ~,
\end{equation}
where $P(z)$ is the power of backscattered light from altitude $z$, $d(z)$ is a known function and $\alpha(x)$ is the unknown extinction coefficient, which can be assumed to be non--negative.

In the last decade, several methods have been proposed for the solution of this problem, including Tikhonov regularization \cite{shcherbakov2007regularized, pornsawad2008ill}, Levenberg--Marquardt \cite{pornsawad2012retrieval} and Richardson--Lucy \cite{garbarino2016expectation}.
These approaches have represented a substantial improvement with respect to previously used restoration techniques, which consisted in computing the numerical derivative of the smoothed data \cite{klett1981stable}:
first, because they explicitly cast the problem as an inverse problem and use explicit regularization (Tikhonov) or an early stopping technique (Levenberg--Marquardt and Richardson--Lucy) to avoid excessive backprojection of noise on the estimated extinction profile;
second -- this only holds true for Richardson--Lucy, although it could be easily introduced in Levenberg--Marquardt -- because they impose a non--negativity constraint, which further helps to reduce instability. On the other hand, to the best of our knowledge, there has been no attempt so far to properly account for noise statistics in the inversion procedure. Indeed, all the proposed methods work on the log--transformed problem
\begin{equation}
y(z) := \log \left( \frac{d(z)}{P(z)} \right) = \int_0^z \alpha(x) dx
\end{equation}
and noise on $y(z)$ is treated either as Gaussian with constant variance (in plain Tikhonov and Levenberg--Marquardt), or as Poisson (in Richardson--Lucy). However, neither of these is a realistic assumption: in general, it is reasonable to assume that $P(z)$ is Poisson, and therefore $y(z)$ is not; neither it is reasonable to assume that the variance of $y(z)$ is the same at all altitudes.

In order to account for the noise statistics, in this study we solve the inverse problem in its non--linear form.
We consider two frameworks, one being early stopped maximum likelihood, the other one with an explicit $\ell^2$ regularization term. 
Specifically, we derive two iterative algorithms by using the Karush-Kuhn-Tucker conditions \cite{zangwill1969nonlinear}, so as to include a non--negativity constraint for the solution.
We use synthetic and experimental data to show that the extinction profiles recovered by the proposed algorithms are substantially better than those recovered by state--of--the--art methods.
In order to further investigate the effects of modelling noise statistics,  we also apply a weighted Tikhonov algorithm to the log--transformed problem and show that this small modification (as opposed to plain Tikhonov) does reduce the instability, even though is not as effective as the newly proposed algorithms. 

The plan of the paper is as follows: in Section 2 we provide a brief introduction to the lidar inverse problem; in Section 3 the proposed algorithms are derived and their properties are described; in Section 4 we present the numerical experiments on synthetic data and in Section 5 an example of application to experimental data. Our conclusions are offered in Section 6.

\section{The Raman lidar equation}

In the atmospheric application of lidar technology, a laser emitted at wavelength $\lambda$ is headed towards the sky; when the backscattered light is measured at the same wavelength, as a consequence of \textit{elastic} scattering, the signal $P_{\lambda}(z)$ from altitude $z$ can be expressed as
\begin{equation}\label{a1_E}
 P_{\lambda}(z) = \frac{C_{\lambda}}{z^2} \beta_{\lambda}(z) e^{\displaystyle -2 \int_0^z \alpha_{\lambda}(x) dx}~.
\end{equation}

Here $C_{\lambda}$ is a known constant accounting for the laser optical power and the overall detection efficiency, the term $\dsp \frac{1}{z^2}$ is the signature of standard signal decay from the backscattering altitude $z$ and $\beta_{\lambda}(z)$ and $\alpha_{\lambda}(z)$ are the (non--negative) backscattering and extinction coefficient, respectively. In this equation both $\alpha_{\lambda}(z)$ and $\beta_{\lambda}(z)$ have to be regarded as unknown. However, the backscattered light can also be measured at a different wavelength $\mu$, as a consequence of inelastic (Raman) scattering. In this case, the backscattering can be entirely addressed to the \textit{molecular} contribution \cite{ansmann1990measurement}, which can be assumed to be proportional to the molecular density profile of the atmosphere $\rho(z)$, and the signal takes the form
\begin{equation}\label{a1}
P_{\mu|\lambda}(z) = \frac{C_{\mu}}{z^2} \rho(z) e^{\displaystyle - \int_0^z \bigl( \alpha_{\mu}(x) + \alpha_{\lambda}(x) \bigr)dx }~,
\end{equation}
where $C_{\mu}$ includes laser optical power, detection efficiency and Raman scattering cross--section and $\alpha_{\mu}(z)$ is the extinction coefficient at the backscatter wavelength $\mu$. The advantage of Raman equation is evident: since $\rho(z)$ can be determined from the Rayleigh theory of scattering,  the total extinction coefficient $\alpha(z) := \alpha_{\lambda}(z) + \alpha_{\mu}(z)$ can be estimated.
Furthermore, the relation between the particles' contribution at the two different wavelengths can be derived from Mie scattering theory \cite{ansmann1992independent}, allowing one to disentangle the two different contributions, just by their sum. 
Most modern lidar systems have both elastic and vibrational Raman channels; therefore, in this study we focus on estimating $\alpha(z)$ from Raman measurements 
(the estimate of $\beta_{\lambda}(z)$ from elastic measurements, once $\alpha_{\lambda}(z)$ is known, is straightforward).

As a first step, we conveniently rewrite equation (\ref{a1}) as 
\begin{equation}
P (z)  = d(z) e^
{\displaystyle -
\int_0^z \alpha(x)\,  dx
}~,
\label{A3ZZ}
\end{equation} 
where subscripts ${\lambda,\mu}$ are omitted but implied, and 
the function $d(z)$ is defined as 
\begin{equation}
d (z)  := \frac{C_{\mu} }{z^2} \rho(z)> 0.
\label{d_definition}
\end{equation} 

We now discretize equation (\ref{A3ZZ}) by introducing two (a priori independent) sets of points $\{z_i\}_{i=1,\dots,N}$ and $\{x_j\}_{j=1,\dots,M}$, and defining:
\begin{itemize}
	\item the noise--free measurement vector $\textbf{P}^*$ such that $P_i^* := P(z_i)$, $i=1,\dots,N$;
	\item the known system function vector $\textbf{d}$ such that $d_i := d(z_i)$, $i=1,\dots,N$;
	\item the unknown exact solution vector $\bm{\alpha}^*$ such that $\alpha_j^* := \alpha(x_j)$, $j=1,\dots,M$;
	\item the linear operator $\textbf{L}$, discretizing the integral.
\end{itemize}
Therefore, we rewrite (\ref{A3ZZ}) in the following discrete compact form 
\begin{equation}
\textbf{P}^*  = \textbf{d} \odot e^{ - \textbf{L} \bm{\alpha}^* } ~,
\label{signal_exact}
\end{equation} 
where $\odot$ is the Hadamard product.
From here on, we assume that the noisy data $\textbf{P}$ have a Poisson distribution, whose mean value is given by the noise--free signal $\textbf{P}^*$ (\ref{signal_exact}). In addition, we assume that noise values at different altitudes are independent, so that the likelihood of a given (permissible) extinction profile $\bm{\alpha}$ is given by
\begin{equation}
 p(\textbf{P}|\bm{\alpha}) = \prod_{i=1}^N p (P_i|\bm{\alpha}) = \prod_{i=1}^N  e^ { \displaystyle - \left( d_i \; e^{- (\textbf{L} \bm{\alpha})_i }\right) } \frac{ \left(d_i \; e^{ - (\textbf{L} \bm{\alpha})_i} \right)^{P_i }}{P_i !}~
\label{eq:Poisson_1}
\end{equation} 
and its domain of definition, i.e. the set of the permissible objects, is the non--negative orthant.
Therefore, our inverse problem is the one of retrieving a non--negative estimate of the exact solution $\bm{\alpha}^*$, given the recorded signal $\textbf{P}$.

\section{Iterative algorithms for the inverse problem}

In this Section we introduce two inversion algorithms to solve the lidar inverse problem, by taking into account the Poisson statistics of the recorded signal $\textbf{P}$. 
The first one belongs to the class of the early stopped maximum likelihood methods; in this context, it is often more convenient to maximize the log--likelihood instead of the likelihood itself, i.e.:
\begin{equation}
\dsp l({\bm{\alpha}}) := \log(p(\textbf{P}|\bm{\alpha}))= \sum_{i=1}^N {\Bigl( P_i \log(d_i) - (\textbf{L}\bm{\alpha})_iP_i  - d_i e^{ - (\textbf{L}\bm{\alpha})_i} - \log (P_i !) \Bigr)},
\label{A7}
\end{equation} 
where the algorithm is stopped at an appropriate iteration, before reaching convergence, in order to prevent excessive backprojection of noise onto the reconstruction.
A second strategy to solve the inverse problem is to supplement the log--likelihood function with a proper penalty term (\textit{prior} in Bayesian terms), accounting for the regularization; the most natural choice in our application is to add an $\ell^2$ penalty term weighed by a regularization parameter $\gamma>0$, and then devise an algorithm to maximize the penalized log--likelihood 
\begin{equation}
S(\bm{\alpha}) := l(\bm{\alpha}) - \gamma \| \bm{\alpha} \|_2^2~.
\label{pen_like}
\end{equation}

In the following, we use the Karush-Kuhn-Tucker (KKT) conditions to derive two iterative algorithms that maximize the log--likelihood (\ref{A7}) and the penalized log--likelihood (\ref{pen_like}), respectively.

\subsection{KKT: Maximum Likelihood for Poisson data}
The first iterative algorithm we propose aims at maximizing (\ref{A7}). Before proceeding with the derivation of the algorithm, we first prove that function (\ref{A7}) has a unique maximum in its domain of definition.\\

\textit{Existence and uniqueness of the maximum}. The maximum of the log--likelihood exists and is unique, provided that $N\ge M$, $\textbf{L}$ is full rank and $P_N>0$.
By elementary computations we obtain that the gradient and the Hessian matrix of $l({\bm{\alpha}})$ are given by 
\begin{equation}
({\bf {\nabla}}l )_j ( {\bm{\alpha}}) = \frac{\partial l ({\bm{\alpha}})}{\partial \alpha_j} = \sum_{i=1}^N {\Bigl(  L_{i,j} { \bigl( d_i  e^{ - (\textbf{L}\bm{\alpha})_i}} - P_i \bigr) \Bigr)}~,
\label{A10}
\end{equation}
and
\begin{equation}
\bigl[ \textbf{H} ({\bm{\alpha}}) \bigr]_{k,j} = \frac{\partial^2 l ({\bm{\alpha}})}{\partial \alpha_k \partial \alpha_j} = - \sum_{i=1}^N { \Bigl( L_{i,j}  L_{i,k} { d_i  e^{ - (\textbf{L}\bm{\alpha})_i}} \Bigr)}~,
\label{A8}
\end{equation} 
for $j, k=1, ..., M$.

From equation (\ref{A8}), it follows that, for any vector $\bm{\mathit{v}} \in \mathbb{R}^N$  the relation
\begin{equation}
\bm{\mathit{v}}^T \textbf{H} ({\bm{\alpha}})  \bm{\mathit{v}} 
= \sum_{k,j=1}^N { \Bigl( v_k \bigl[ \textbf{H} ({\bm{\alpha}}) \bigr]_{k,j} v_j \Bigr)} 
= - \sum_{i=1}^N {\Bigl( d_i  e^{ - (\textbf{L}\bm{\alpha})_i}  (\textbf{L} \bm{\mathit{v}} )_i^2 \Bigr)}
\label{A9bis}
\end{equation} 
holds. Since $d_i > 0$ $\forall$ $i = 1,..,N$, assuming that $N\ge M$ and $\textbf{L}$ is full rank, we have
$\bm{\mathit{v}}^T \textbf{H} ({\bm{\alpha}})  \bm{\mathit{v}} < 0 \;\; \forall \bm{\mathit{v}} \ne \bm{0}$;
therefore the Hessian matrix is negative definite, implying $l (\bm{\alpha})$ is strictly concave.
We now observe that the likelihood can be written in the following way, up to constant terms:
\begin{eqnarray}\label{eq:15}
l({\bm{\alpha}}) & \sim  - \sum_{i=1}^N { \Bigr( (\textbf{L}\bm{\alpha})_iP_i  + d_i e^{ - (\textbf{L}\bm{\alpha})_i} \Bigr)} \nonumber\\ 
& = - \left(\sum_{j=1}^M {\Bigr( (\textbf{L}^T \textbf{P})_j \alpha_j \Bigr)}+\sum_{i=1}^N {\Bigr( d_i e^{ - (\textbf{L}\bm{\alpha})_i} \Bigr)}\right)
\end{eqnarray} 
and it is defined on the non--negative orthant ($\bm{\alpha} \ge \bm{0}$).

From here on, we assume that the last component $P_N$ of the data $\textbf{P}$ is non--zero (otherwise, it is sufficient to cut the data at the last non--zero component).
In addition, $\textbf{L}$ discretizes the cumulative integral from the first height, therefore $\textbf{L}^T$ can be thought of as a cumulative, weighted sum, from the last height. These two facts guarantee that $( \textbf{L}^T \textbf{P} )_j$ never vanishes (it is actually always positive).
We now observe that
\begin{equation}\label{eq:3rispref}
0<\sum_{i=1}^N {\left( d_i e^{ - (\textbf{L}\bm{\alpha})_i} \right)} \le \sum_{i=1}^N d_i  ~;
\end{equation}
moreover, since $(\textbf{L}^T \textbf{P})_j \ne 0$ $\forall$
$j = 1,..,M$, if we take a vector $\bm{\alpha}$ such that $\alpha_j\to+\infty$, for any $j$, we have $\sum_{j=1}^M {\left( (\textbf{L}^T \textbf{P})_j \alpha_j \right)} \to +\infty$.
Therefore, thanks to the first inequality in (\ref{eq:3rispref}), we have
\begin{equation} 
\lim_{\|{\bm{\alpha}\| \to +\infty}}  l({\bm{\alpha}}) = -\infty.
\label{coerc}
\end{equation} 
This fact, together with the strict concavity, guarantees that the likelihood admits a unique global maximum point (see Lemma 2 in \cite{lange1995globally}).\\

\textit{Derivation of the algorithm}. In general, the KKT conditions are first--order necessary conditions for optimality under given constraints. In our context,
for the maximum $\tilde{\bm{\alpha}}$ of $l (\bm{\alpha})$, they can be expressed as
\begin{equation}
\cases{
\alpha_j \frac{\partial l (\bm{\alpha})}{\partial \alpha_j} \Big|_ {\bm{\alpha}=\tilde{\bm{\alpha}}}^{} = 0 & \text{$j = 1,\dots, M$} \\ 
\frac{\partial l (\bm{\alpha})}{\partial \alpha_j} \Big|_ {\bm{\alpha}=\tilde{\bm{\alpha}}}^{} \le 0 & \text{if $\tilde{\alpha}_j=0$}~.
}
\label{A11}
\end{equation}
Importantly, since $l (\bm{\alpha})$ is concave, in our case the KKT conditions are not only necessary but also sufficient conditions for an object $\tilde{\bm{\alpha}}$ to be 
the maximum of $l (\bm{\alpha})$.

By combining equation (\ref{A10}) with equations (\ref{A11}) we 
obtain that the maximum of $l (\bm{\alpha})$ must 
be solution of the following non--linear equation 
\begin{equation}
\tilde{\alpha}_j = \frac{ \bigl( \textbf{L}^T ( \textbf{d }\odot  e^{ - \textbf{L} \tilde{\bm{\alpha}}} ) \bigr) _j  }{( \textbf{L}^T \textbf{P} )_j}     \tilde{\alpha}_j \quad j = 1,\dots, M~,
\label{A14}
\end{equation} 
where, as pointed out above, the denominator is always positive.
Let us now introduce
the non--linear operator 
$\textbf{T}: \mathbb{R}^N \to \mathbb{R}^N$
\begin{equation}
\textbf{T} (\bm{\alpha}) := \frac{\textbf{L}^T ( \textbf{d } \odot e^{ - \textbf{L} \bm{\alpha} } )}{\textbf{L}^T \textbf{P}} \odot\bm{\alpha}~,
\label{A15}
\end{equation} 
where the division has to be intended as the inverse of the Hadamard product, replaced with the usual division symbol for the sake of readability.
Therefore, it follows that the maximum of $l (\bm{\alpha})$ must 
be fixed point of $\textbf{T}$, i.e.
\begin{equation}
\tilde{\bm{\alpha}} = \textbf{T}(\tilde{\bm{\alpha}})~.
\end{equation}
By applying the method of successive approximations (as done, for instance, in \cite{lanteri2002penalized}), to this fixed point equation, we obtain the following iterative algorithm for the maximization of $l(\bm{\alpha})$
\begin{equation}
{\bm{\alpha}}^{( n+1 )} :=  \textbf{T} ({\bm{\alpha}}^{( n )}) = \frac{ \textbf{L}^T ( \textbf{d } \odot e^{ - \textbf{L} {\bm{\alpha}}^{( n )} } )  }{ \textbf{L}^T \textbf{P} } \odot{\bm{\alpha}}^{( n )}~.
\label{kkt}
\end{equation} 

We notice that there are no obvious guarantees that algorithm (\ref{kkt}) will converge; however, it is possible to interpret (\ref{kkt}) as follows
\begin{equation}
{\bm{\alpha}}^{( n+1 )} =  {\bm{\alpha}}^{(n)} + \bm{D}({\bm{\alpha}}^{(n)}) \nabla l({\bm{\alpha}}^{(n)}) ~.
\label{kkt_gradient}
\end{equation}
Now this is a scaled gradient algorithm with a constant and unitary step size, and a diagonal scaling matrix $\bm{D}({\bm{\alpha}}) :=  diag \Bigl( \displaystyle \frac{\bm{\alpha}}{\textbf{L}^T \textbf{P}} \Bigr)~$.
For this class of algorithms convergence is guaranteed by introducing a line search procedure \cite{lanteri2001general,bonettini2008scaled}; hence, in the numerical simulations below we have utilized the Armijo rule \cite{armijo1966minimization}, designed to respect the positivity constraint. \\
Even though the numerical results discussed below have been obtained with the scaled gradient algorithm (\ref{kkt_gradient}), we provide here some additional comments on the multiplicative algorithm (\ref{kkt}).

\textbf{Remark 1}. If  all the components of the initial guess $\bm{\alpha}^{( 0 )}$ are (strictly) positive, then such property is preserved by each $\bm{\alpha}^{( n )}$, $n>0$. As a consequence, $\bm{\alpha}^{( n )}$ is a fixed point of the operator $\textbf{T}$ only if it is the maximum of the log--likelihood function. 


\textbf{Remark 2}. If the iterations of the multiplicative algorithm (\ref{kkt}) converge, they converge to the maximum (i.e. the limit point satisfies both the KKT conditions (\ref{A11})); this is a consequence of the concavity of $l(\alpha)$. In fact, in all the numerical experiments performed for this study we did observe numerical convergence for the multiplicative algorithm (\ref{kkt}); these results are not shown in this article.

\textbf{Remark 3}. Algorithm (\ref{kkt}) has been described in \cite{lange1995globally} for an imaging problem in Computerized Tomography, $\textbf{L}$ being the Radon Transform rather than the cumulative integral. The authors discuss the convergence properties of the algorithm. By assuming the existence and the uniqueness of an optimal point  $\tilde{\bm{\alpha}}$ satisfying the strict constraint $\tilde{\alpha}_j > 0$ $\forall j=1,\dots M$, they show that local convergence is guaranteed if $(\textbf{L} \tilde{\bm{\alpha}})_i < 1, \; \forall i=1,\dots N$.

As a final comment, we recall that in this study the iterations will not be performed until convergence. Indeed, due to the ill-posedness of the problem, the maximum likelihood estimate would not provide a sensible reconstruction of the unknown. Even though the theoretical properties of early stopping are still matter of debate, heuristics suggest that early stopping can have a regularizing effect on the solution, and our numerical simulations below confirm this conjecture.

\subsection{KKT\_L2: Penalized maximum likelihood for Poisson data}

We now consider the maximization of the penalized log--likelihood (\ref{pen_like}).
Since we adopt the very same strategy used for the maximization of (\ref{A7}) we are not going to repeat the derivation here. In fact, thanks to the regularizing term, function (\ref{pen_like}) is strictly concave even without the requirements $N\ge M$ and $\textbf{L}$ full rank; in addition, it has the same limiting behavior (\ref{coerc}) of  the log--likelihood $l(\bm{\alpha})$; consequently, it admits a unique global maximum in its domain.
 By writing down the KKT conditions with the non--negativity constraint, a similar iterative algorithm is derived, which reads as follows:
\begin{equation}
{\bm{\alpha}}^{( n+1 )} =  \frac{ \textbf{L}^T ( \textbf{d } \odot e^{ - \textbf{L} {\bm{\alpha}}^{( n )} } )  }{ \textbf{L}^T \textbf{P} + 2 \gamma \bm{\alpha}^{(n)}} \odot{\bm{\alpha}}^{( n )}~,
\label{kkt_l2}
\end{equation} 
where $\gamma$ is the regularization parameter. Again, we interpret the algorithm as a scaled gradient one and apply the line search procedure to ensure convergence. 
We notice that the considerations expressed in Remarks 1 and 2 hold also for the multiplicative version of this algorithm. In addition, we remark the following facts.

\textbf{Remark 4}. Since regularization is obtained here by introduction of a regularizing term in the functional, there is no need to stop the iterations in advance. In the numerical sections, we will perform 200 iterations, and always observe convergence of the algorithm.

\textbf{Remark 5}. In standard Tikhonov regularization, the regularization parameter can be interpreted as the ratio between the variance of the prior and the variance of the likelihood function. 
In our model, since the likelihood is Poisson, the regularization parameter is interpretable as the inverse variance of the prior distribution only.
Therefore the optimal regularization parameter should not depend on the signal--to--noise ratio. 
\section{Numerical experiments}

\begin{figure}
\begin{center}
\includegraphics[width=10cm]{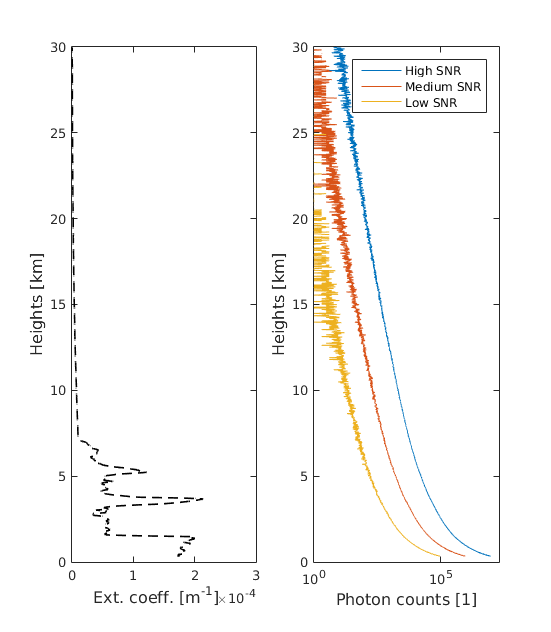}
\caption{Left panel: the true extinction profile for Simulation 1. Right panel: synthetic data with Poisson noise at three different SNRs (logarithmic scale on the $x$ axis).}
\label{fig:truedata}
\end{center}
\end{figure}

\begin{figure}
\begin{center}
	\includegraphics[width=13.4cm]{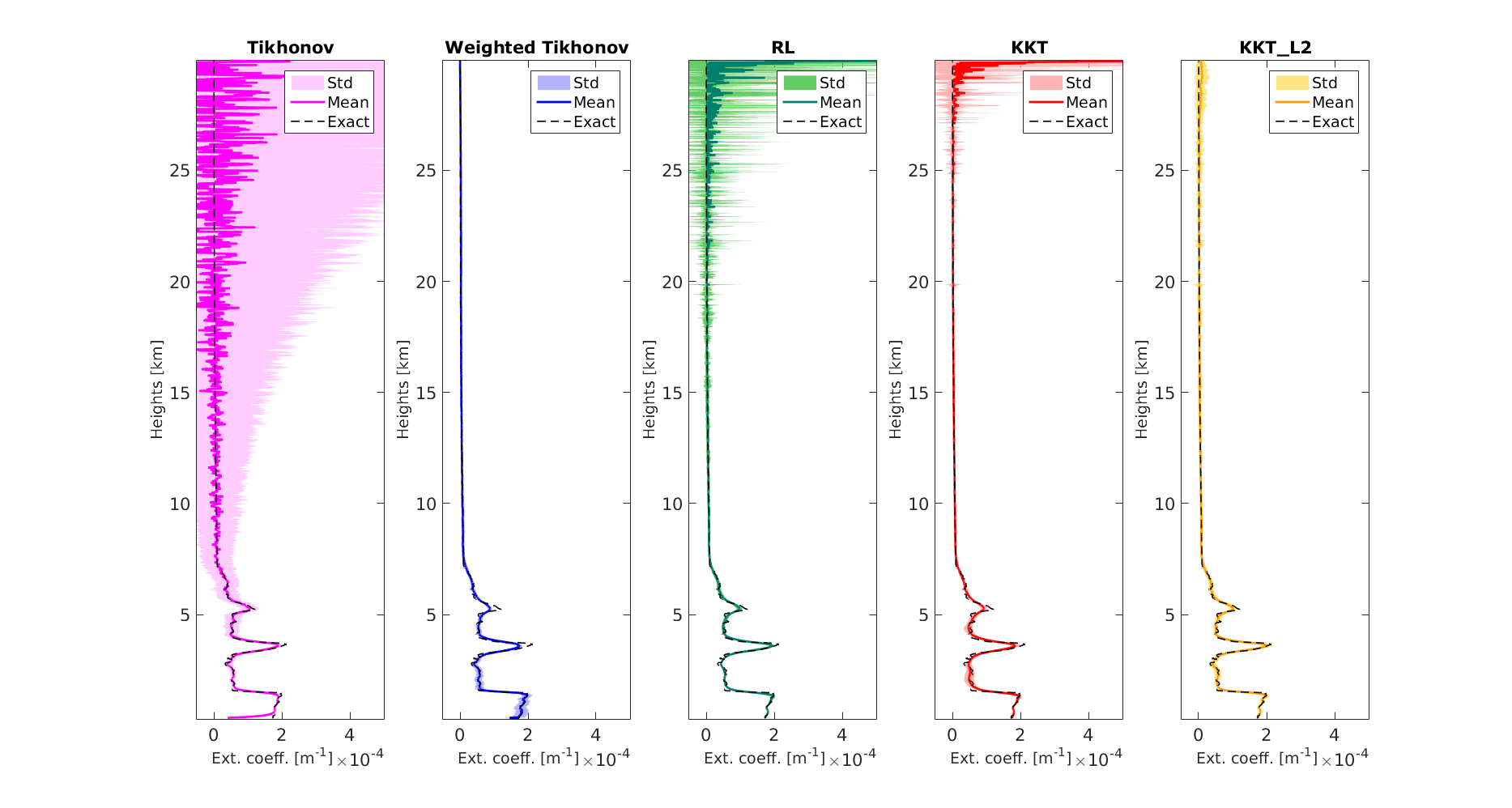}
	\includegraphics[width=13.4cm]{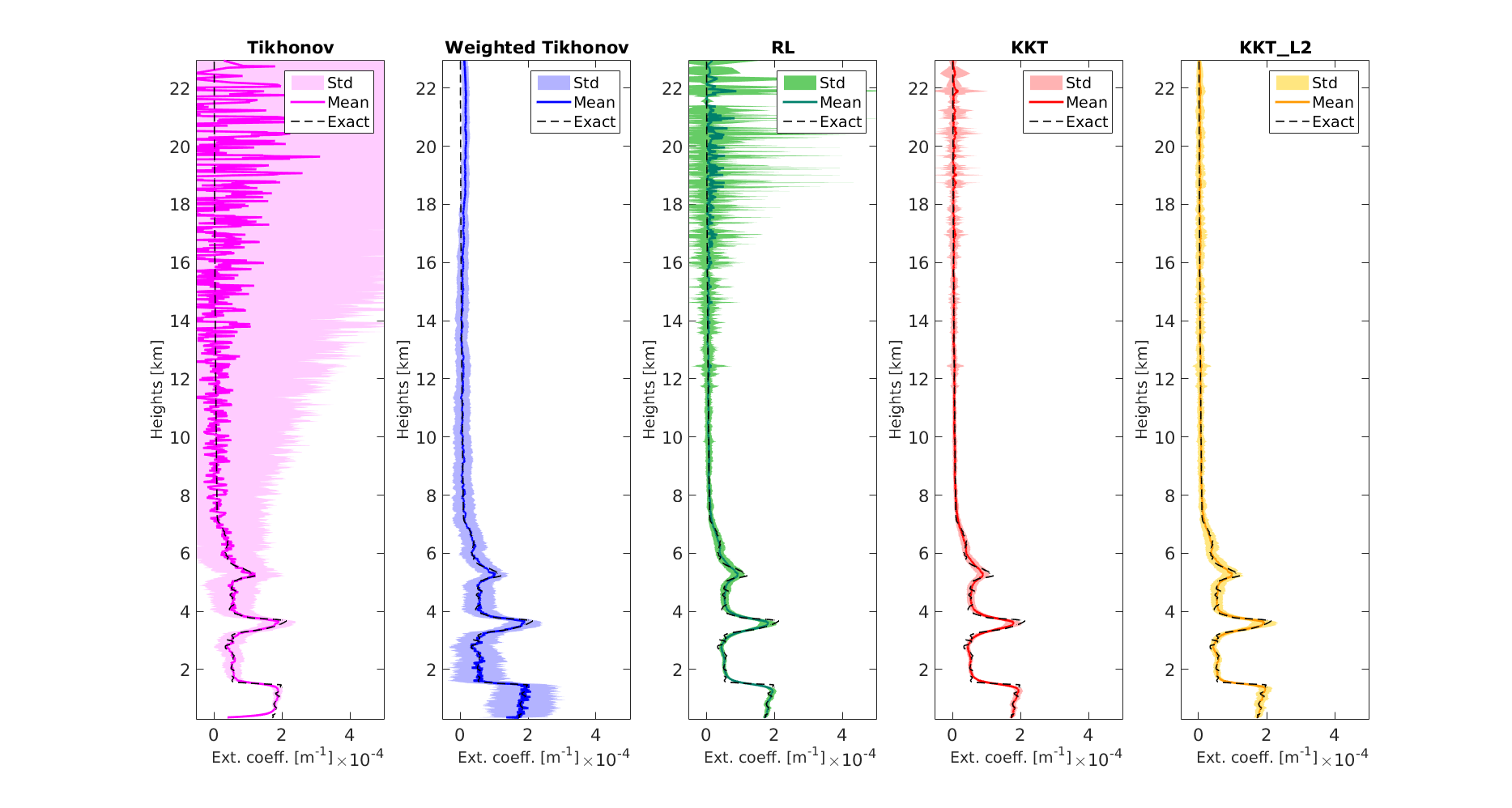}
	\includegraphics[width=13.4cm]{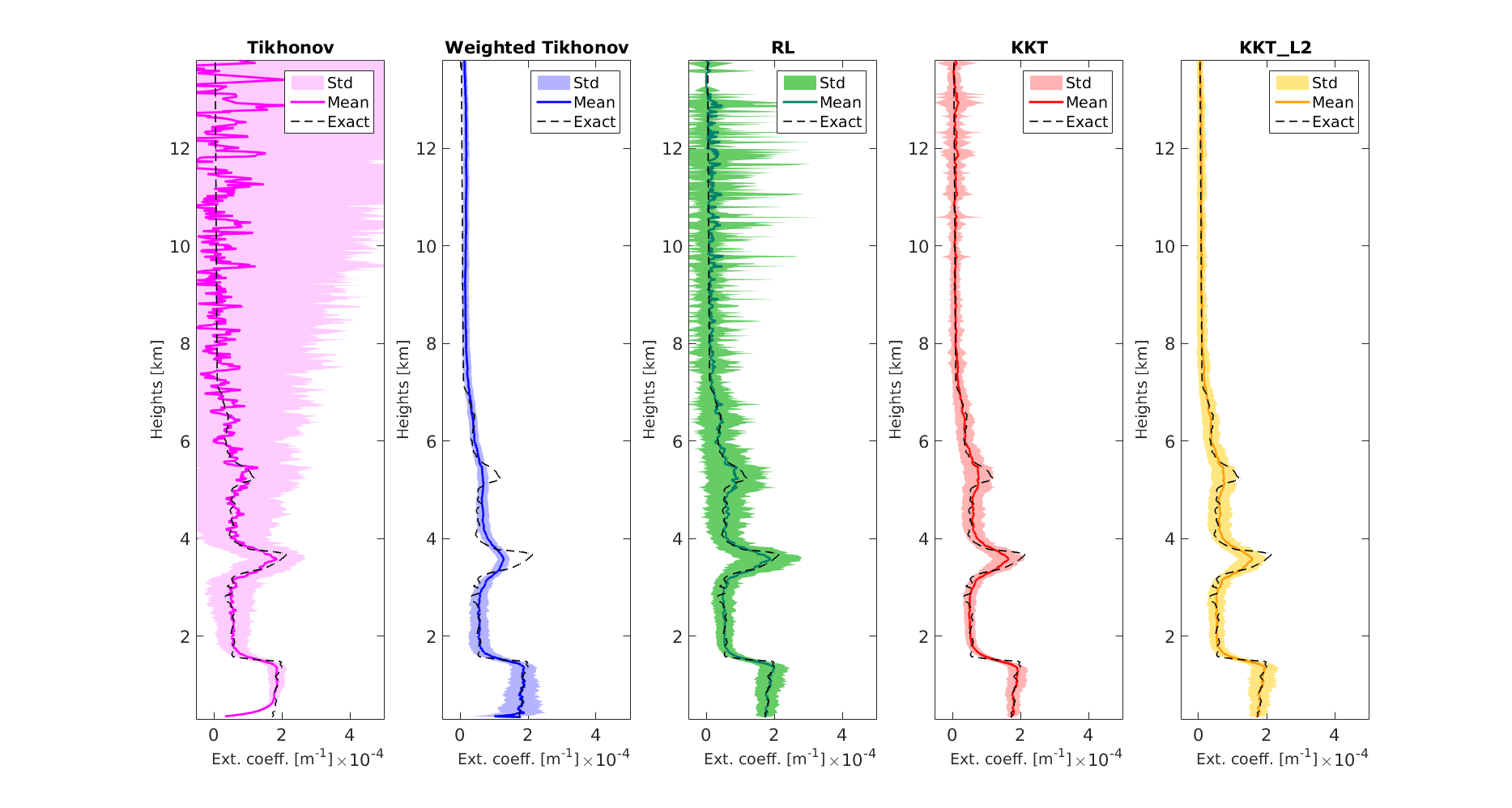}
	\caption{Simulation 1. Mean and standard deviation of the estimated profiles over 100 noise realizations for high (top), medium (centre) and low (bottom) SNR. The parameters are as follows:
 $\gamma_{Tikhonov} = 5 * 10^{3}$, $\gamma_{w. Tikhonov} = 2*10^{7}$, $n_{RL} = 6000$, $n_{KKT} = 100$, $\gamma_{KKT_2} = 2*10^{6}$ (top); $\gamma_{Tikhonov} = 5*10^{3}$, $\gamma_{w. Tikhonov} = 2*10^{7}$, $n_{RL} = 4000$, $n_{KKT} = 120$, $\gamma_{KKT_2} = 2*10^{6}$ (centre); $\gamma_{Tikhonov} = 10^{4}$, $\gamma_{w. Tikhonov} = 3*10^{7}$, $n_{RL} = 3500$, $n_{KKT} = 140$, $\gamma_{KKT_2} = 4*10^{6}$ (bottom).}
	\label{fig:sim1}
	\end{center}
\end{figure}

\begin{figure}
\begin{center}
	\includegraphics[width=14cm]{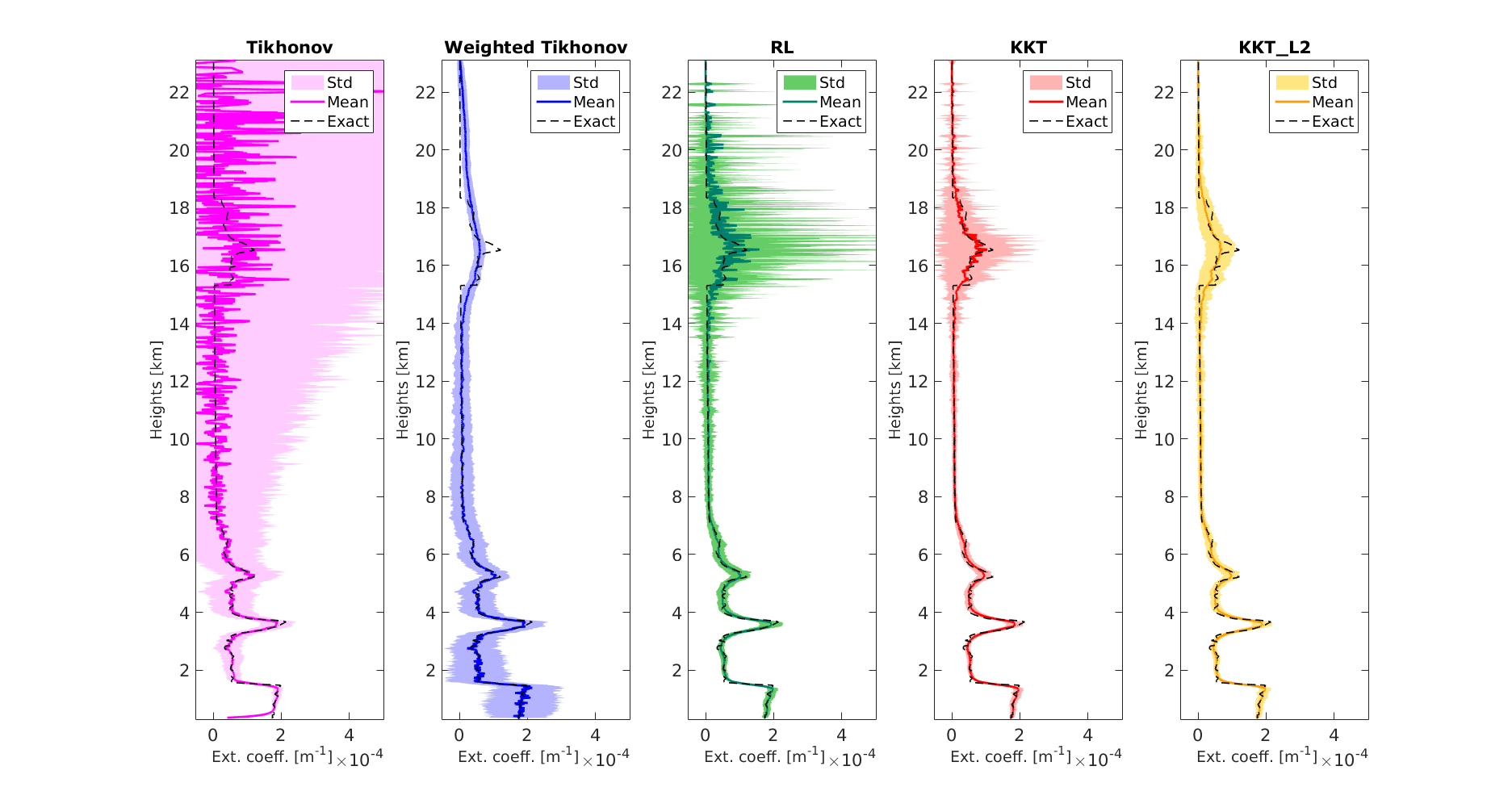}
	\caption{Simulation 2. Mean and standard deviation of the estimated profiles over 100 realizations for medium SNR. Results are obtained with the following parameters:  $\gamma_{Tikhonov} = 10^{5}$, 
	$\gamma_{w. Tikhonov} = 10^{7}$, $n_{RL} = 6000$, $n_{KKT} = 170$, $\gamma_{KKT_2} = 2*10^{6}$.}
	\label{fig:sim2}
	\end{center}
\end{figure}

In this Section we evaluate the proposed algorithms against two sets of synthetic data and compare them to two state--of--the--art methods in this field, namely plain Tikhonov \cite{shcherbakov2007regularized,pornsawad2012retrieval} and Richardson--Lucy  (RL) \cite{garbarino2016expectation}. We omit direct comparison against Levenberg--Marquardt (LM), which was shown in \cite{garbarino2016expectation} to be outperformed by RL, in the same experimental set--up used in this paper (i.e. Raman lidar data).
Since plain Tikhonov assumes that the noise variance is the same at every altitude, but we know that this assumption is incorrect for lidar data, we add to the comparison also a weighted Tikhonov algorithm. We recall that Tikhonov and RL work on the log--transformed problem
\begin{equation}\label{eq:logtrasfdiscr}
\textbf{y} = \textbf{L}\bm{\alpha}~,
\end{equation}
where $\textbf{y} := \log(\textbf{d} / \textbf{P})$.
The weighted Tikhonov solution is then given by
\begin{equation}\label{eq:tiksol}
\tilde{\bm{\alpha}} = \left(\textbf{L}^T\textbf{W}\textbf{L} + \gamma\textbf{I}\right)^{-1}\left(\textbf{L}^T\textbf{W}\textbf{y}\right)~,
\end{equation}
where $\gamma>0$ is a regularization parameter and $\textbf{W}$ is a diagonal weight matrix, accounting for altitude--dependent noise variance: the element $W_{ii}$ represents the inverse of the noise variance at altitude $z_i$. In order to obtain an estimate of these values, we adopt the following strategy: we generate 100 Poisson realizations of $\textbf{P}$, compute the corresponding $\textbf{y}$ and then compute the empirical variance of these 100 realizations of $\textbf{y}$.
For the sake of completeness, we also recall that RL, when applied to the log--transformed data as in \cite{garbarino2016expectation}, is the iterative algorithm whose iterations are described by
\begin{equation}\label{eq:rl}
{\bm{\alpha}}^{(n+1)} = \frac{{\bm{\alpha}}^{(n)}}{\textbf{L}^T \textbf{1}}\odot\textbf{L}^T \frac{\textbf{y}}{\textbf{L}{\bm{\alpha}}^{(n)}}~.
\end{equation}

We now proceed with the description of the two simulated experiments. In each of the two simulations, we start from a true extinction profile and compute the corresponding exact signal through eq. (\ref{a1}). To assess the sensitivity to noise, we then generate 100 Poisson realizations with the same mean, we compute the solutions with all the algorithms, and then show mean and standard deviation of the retrieved extinction coefficients. In order to assess the impact of the signal--to--noise ratio (SNR), we repeat this procedure for three different values of the SNR (high, medium, low). Clearly, as the data are Poisson, the only way to change the SNR is by changing the signal strength; this is done by modifying the constant $C_{\mu}$ in (\ref{d_definition}). 

\textit{Simulation 1}. The first simulation we consider is based on the extinction profile previously used to assess the performances of various inversion algorithms \cite{pappalardo2004aerosol} in the framework of the  European Aerosol Research Lidar Network (EARLINET). Here $M=N=2000$. 

In Figure \ref{fig:truedata} we show the true extinction profile and three examples of synthetic data, for three different noise levels. The extinction profile and the data are plotted according to the standard vertical visualization of lidar quantities: altitudes are on the $y$ axis, the signals and the extinction coefficient on the $x$ axis.
In Figure \ref{fig:sim1} we show the average estimated profiles, with standard deviation, for the five methods and  three different values of the SNR. 
In order to point out the differences between the algorithms, the number of iterations and the regularization parameters have been selected as those that minimize the discrepancy between the mean estimated profile and the true profile, at the lower altitudes.

The results indicate that plain Tikhonov (magenta in the Figure) is extremely affected by noise: when the SNR is high, top row, the profile at low altitudes is correctly reconstructed with almost zero variance; however, the variability of the reconstruction at higher altitudes is large, so that the reconstruction from a single realization (as in a real case) is useless above, say, 10 km; for lower values of the SNR, even the low altitudes present high variability. Of course one could regularize more, but then the estimates at the lower altitudes would be oversmoothed.

Weighted Tikhonov (blue in the Figure) is better: in the high SNR case, the reconstruction is very good, with almost no variance; when the SNR goes down, however, more variance appears at the lower altitudes, and the higher altitudes are over--smoothed. This behaviour is opposite with respect to plain Tikhonov, and suggests that either the Gaussian approximation is not good, or the estimated weights are not optimal.

As discussed in \cite{garbarino2016expectation}, Richardson--Lucy (green in the Figure) improves over Tikhonov: in the high and medium SNR case, the estimated profile is affected by noise only in the higher part of the range; for low SNR, the reconstructions have more variability, as expected, but the profile remains discernible. We observe that the standard deviation in the figure exhibits an irregular behaviour: this is due to the presence of large noise--induced spikes in the single reconstructions, which change from one noise realization to the other one.

The KKT algorithm (red in the Figure) does not seem to depart qualitatively from RL, but it is quantitatively better: while retaining approximately the same mean value, it has lower variance, systematically, at all altitudes and for all three SNR values. This is pleasant and expected, and confirms that the KKT algorithm is better suited to treat the Poisson noise on the signal $\textbf{P}$.

Finally, KKT\_L2 appears to exhibit very good robustness to Poisson noise: all the spikes of RL and KKT disappear here. The reconstructions at high SNR have almost negligible variance and almost no noise, and appear to outperform those of weighted Tikhonov. For lower values of the SNR, we witness some more variance and a little oversmoothing at the lower altitudes, which is expected due to the presence of the regularizing term.\\

\textit{Simulation 2}. Simulation 1 clearly shows that the KKT algorithm improves slightly over RL, but KKT\_L2 improves substantially when it comes to the stability of the solution.  However, it may be objected that KKT\_L2 has a considerable smoothing effect at high altitudes, and therefore it may be questioned whether KKT\_L2 would be able to recover something at high altitudes. To answer this question, in Simulation 2 we take the same synthetic extinction profile of Simulation 1 and add a structure in the higher part of the range.

In Figure \ref{fig:sim2} we show the estimated profiles obtained by the five inverse methods, in the case of medium SNR; similarly to the first simulation, the parameters have been chosen in order to minimize the discrepancy between the mean estimated profile and the exact profile, at the lower heights. 
As expected, KKT\_L2 has a smoothing effect at high altitudes, but its reconstructions remain relatively stable. 
As far as RL and KKT are concerned, they may seem to provide better results, because their mean values follow the true profile more accurately. However, this is due to the fact that they are maximum likelihood methods, and therefore more unbiased than KKT\_L2 (not completely unbiased because they are stopped early); but the variance of their reconstructions is so large (for RL especially) that the reconstruction from a single realization is almost useless in the higher part of the range. 

\begin{figure}
\begin{center}
	\includegraphics[width=14cm]{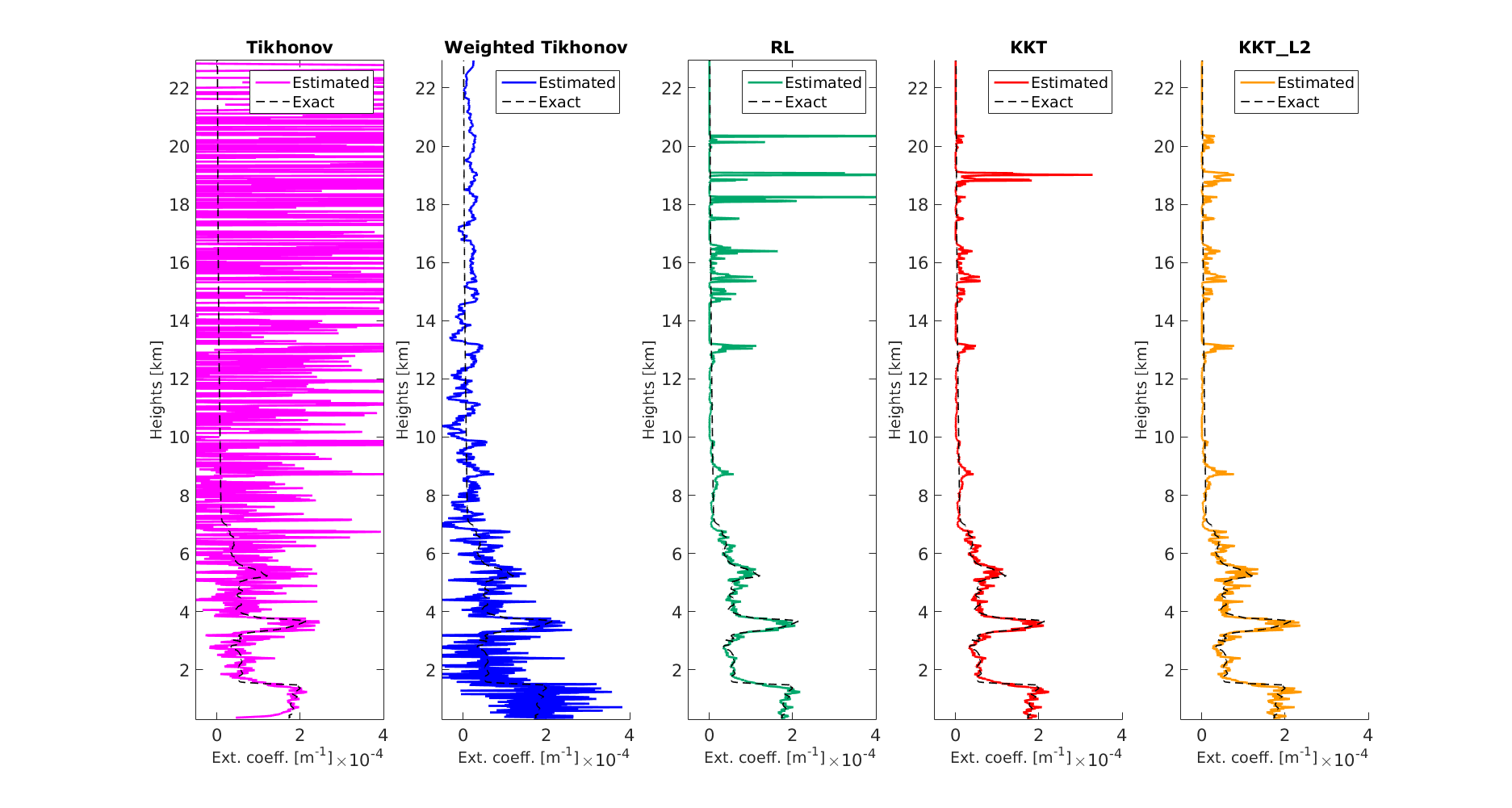}
	\includegraphics[width=14cm]{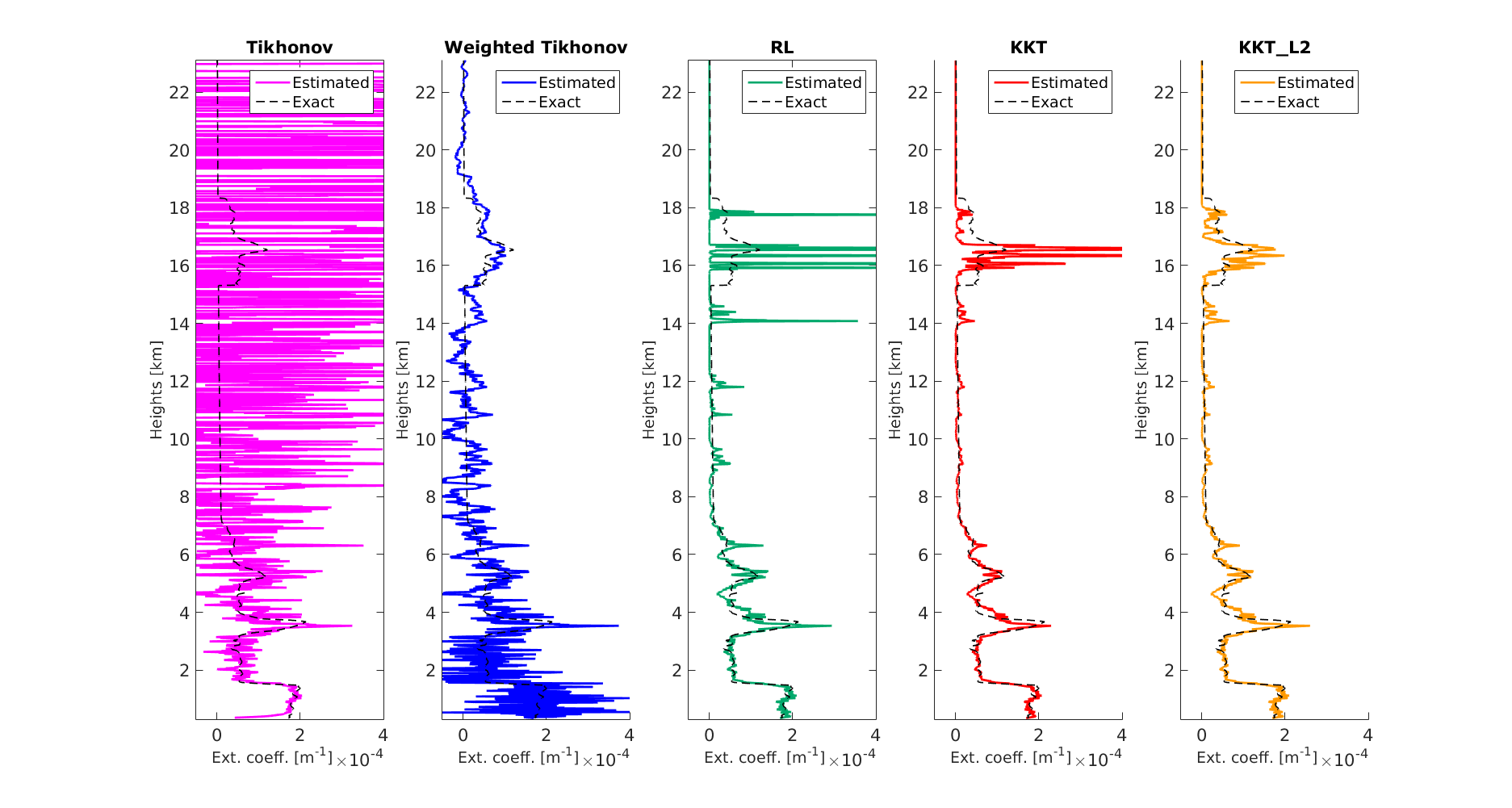}
	\caption{Reconstruction from a single realization for medium SNR: top, Simulation 1; bottom, Simulation 2.}
	\label{fig:sim1_vs_sim2}
	\end{center}
\end{figure}

\begin{figure}
\begin{center}
	\includegraphics[width=14cm]{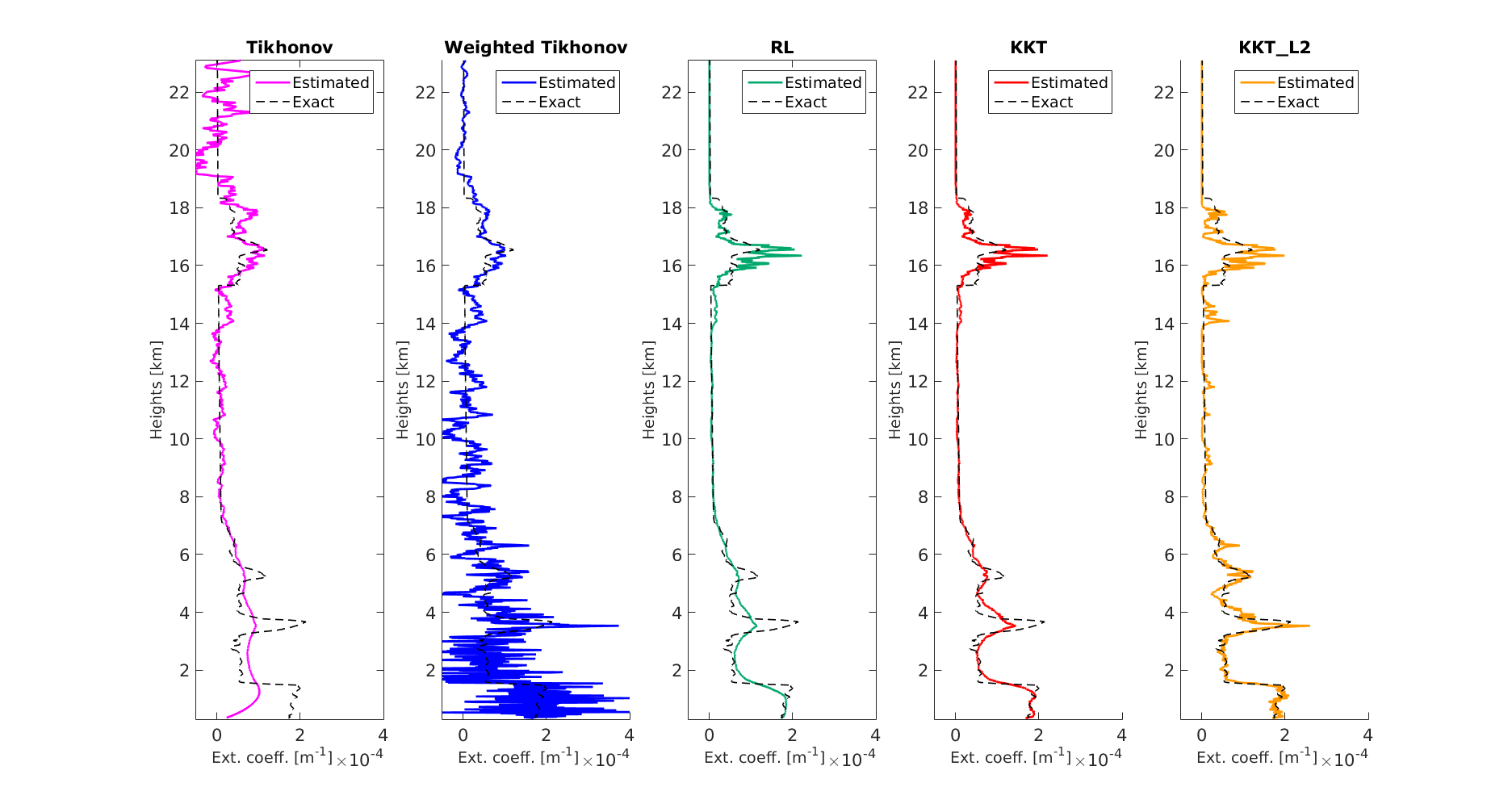}
	\caption{Reconstruction from a single realization for Simulation 2, where the regularization parameters have been set in order to provide a better reconstruction of the structure at high altitudes. Explicitly, the parameters are as follows: $\gamma_{Tikhonov} = 5*10^{5}$,  $\gamma_{w. Tikhonov} = 10^{7}$, $n_{RL} = 350$, $n_{KKT} = 70$, $\gamma_{KKT_2} = 2*10^{6}$. }
	\label{fig:sim2_highaltitudes}
	\end{center}
\end{figure}

To highlight this point, in Figure \ref{fig:sim1_vs_sim2} we present the estimated profiles from a single realization: in the top row we plot the results from Simulation 1, in the bottom row from Simulation 2. Clearly, it would be difficult to tell which is which by looking either at plain Tikhonov, where oscillations completely destroy the information content of the signal, or Richardson--Lucy, where similar spiky structures appear at around the same altitudes. Weighted Tikhonov manages to recognize
the presence of the structure on the top, but the retrieved profile at the lower altitudes is extremely noisy.
KKT is better than RL, as it produces less spikes when the true underlying profile is flat and seems to recognize the presence of the peak, although over--estimating its intensity and, finally, KKT\_L2 manages to provide a substantially different reconstruction, with a clear indication of the presence and of the intensity of a structure in the higher part of the range.

It is worth to observe that a better estimate of the structure at high altitudes can be obtained also by RL and KKT and, to a lesser degree, plain Tikhonov, just modifying the stop iterations and the regularization parameter. As shown in Figure \ref{fig:sim2_highaltitudes}, the profiles estimated at high altitudes by RL and KKT are much closer to the true solution, and very similar to the one reconstructed by KKT\_L2; however, this improvement at high altitudes comes at the price of a strong oversmoothing in the lower part of the range, particularly for RL.

\section{Application to experimental data}
In this section we compare the performances of all the algorithms using a set of experimental data. 

The signal was acquired on June 17th 2013 by the MALIA (Multiwavelength Aerosol LIdar Apparatus) lidar \cite{boselli2009atmospheric} installed in Naples.
The laser wavelength was 355 nm, with a pulse energy of 100 mJ, 20 Hz repetition rate and a 30 cm aperture telescope. The recordings took place during a Saharan dust event, whose traces can be recognized in the estimated profile. Here $M=N=1000$.

\begin{figure}
\begin{center}
	\includegraphics[width=14cm]{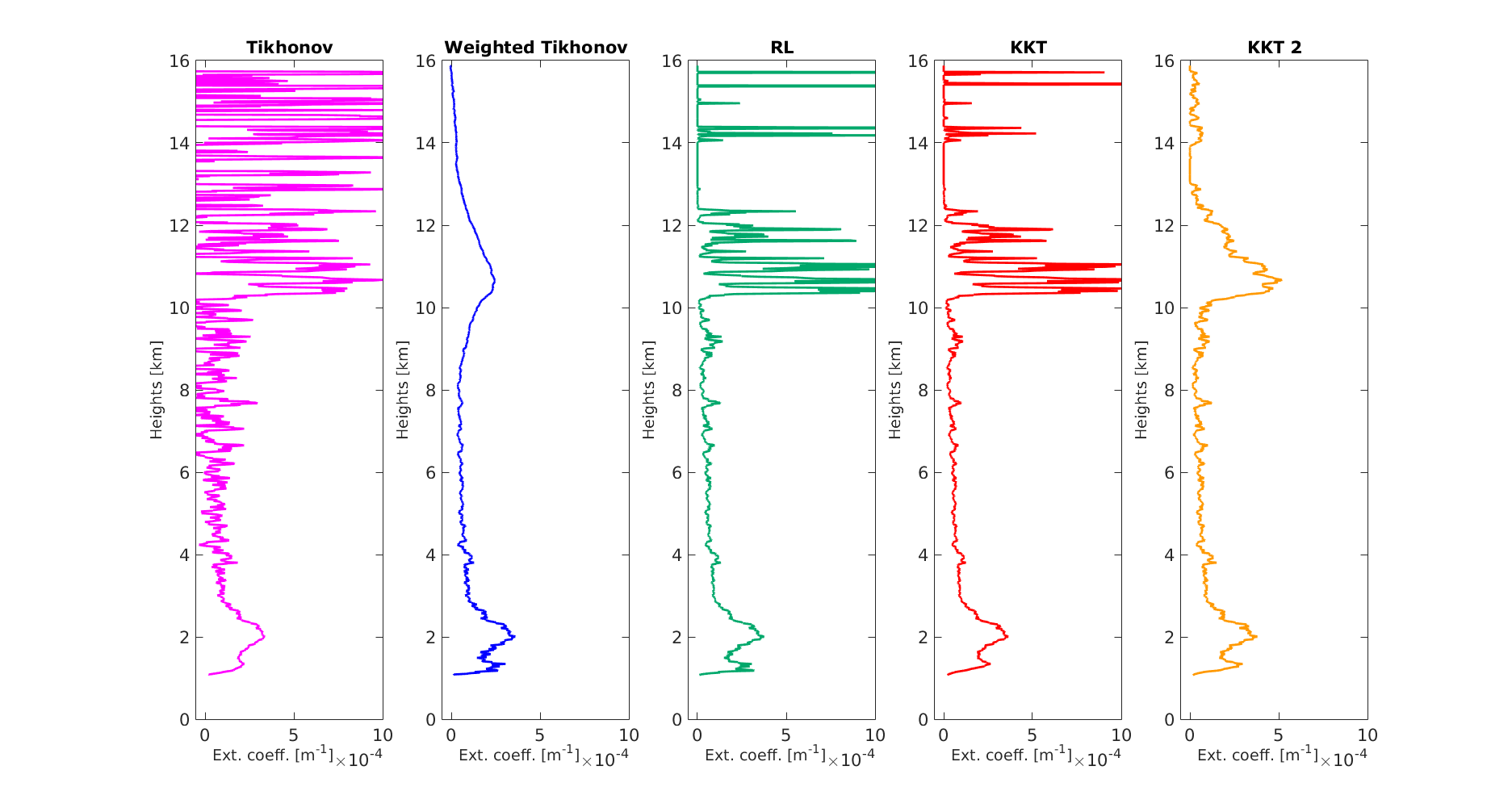}
	\includegraphics[width=14cm]{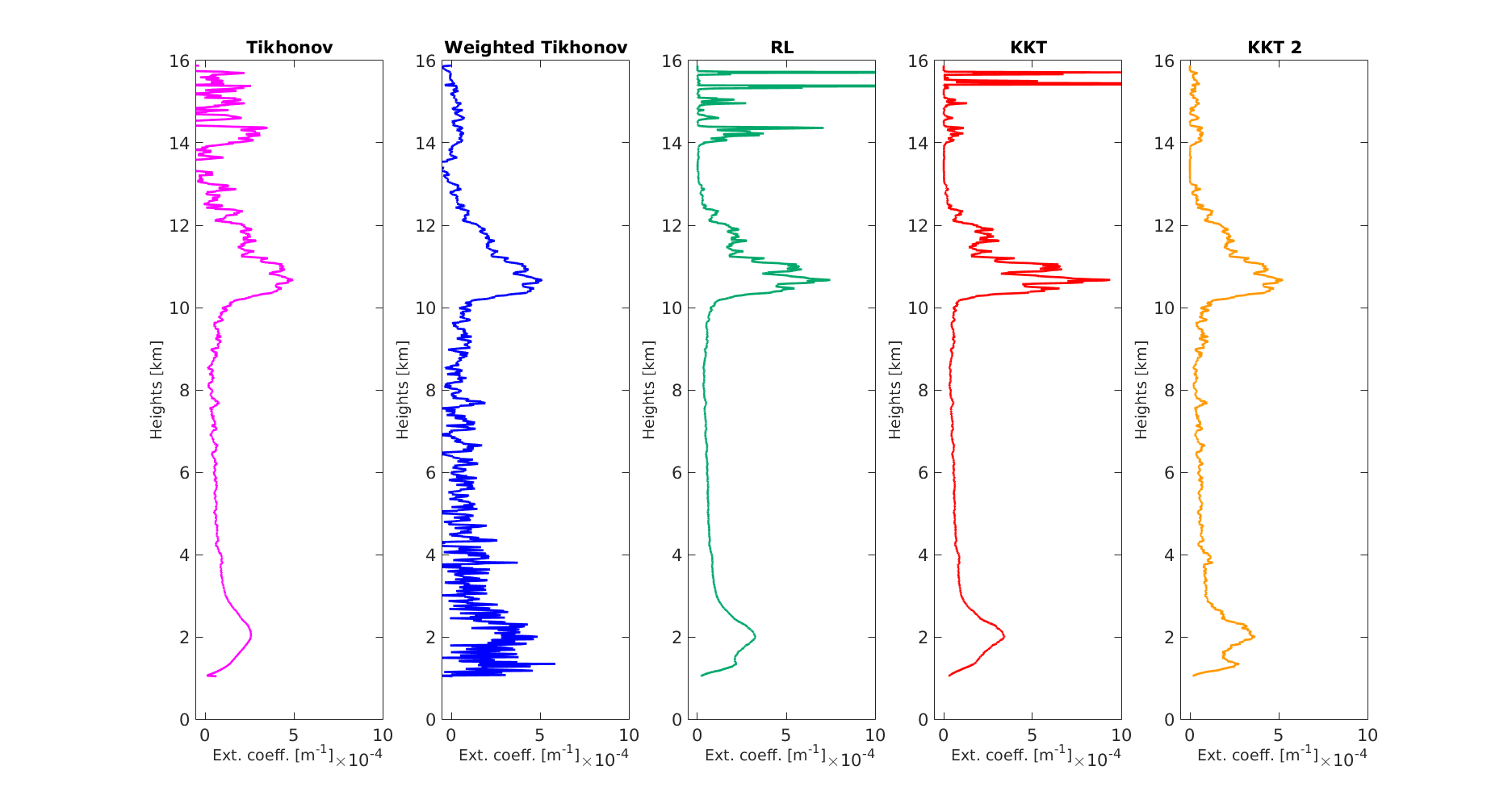}
	\caption{Real data. In the first row the results obtained with the following parameters:  $\gamma_{Tikhonov} = 10^{4}$, $\gamma_{w. Tikhonov} = 10^{8}$, $n_{RL} = 2000$, $n_{KKT} = 170$, $\gamma_{KKT_2} = 5*10^{6}$. In the second row the results obtained with the following parameters:  $\gamma_{Tikhonov} = 10^{5}$, $\gamma_{w. Tikhonov} = 10^{7}$, $n_{RL} = 300$, $n_{KKT} = 40$, $\gamma_{KKT_2} = 5*10^{6}$. We point out that the regularization parameter used for KKT\_L2 is the same in both cases.}	
	\label{fig:real1}
\end{center}
\end{figure}

In Figure \ref{fig:real1} we present the reconstructed extinction coefficients.
Following the results from Simulation 2 above, we present two different estimated profiles, for two different choices of the regularization parameters: one (top row) where the parameters have been selected to obtain similar results (across methods) at low altitudes, one (bottom row) where the parameters have been selected to obtain similar results at high altitudes. 

The results are consistent with those obtained on synthetic data. In the top row, the algorithms obtain very similar reconstructions at the lower altitudes, indicating that the information content of the data is high, and the different priors do not play a major role. In particular, plain Tikhonov provides a good reconstruction at the low heights but wide oscillations at high altitudes; weighted Tikhonov has the opposite behaviour, providing a very smooth reconstruction at high altitudes; the profile estimated by KKT appears to remain reliable up to around 10 kilometres, but presents high peaks from there on; RL is similar to KKT; KKT\_L2 seems to provide a reliable reconstruction across the whole range. 

In the second row, the estimated profiles present similar features but are not identical at high altitudes, indicating that the signal here is less informative, and the prior becomes more relevant.
With this choice, the profiles estimated by plain Tikhonov, RL and KKT are overly smoothed at low altitudes: the general structure is preserved, but the fine details are lost. The profile reconstructed by weighted Tikhonov has, again, the opposite behaviour, becoming too oscillating at low altitudes. 
Again KKT\_L2 provides the best solution: a good balance between fine resolution and smoothness is achieved uniformly at all the altitudes.

\section{Discussion and conclusion}

In this article we have introduced two methods, KKT and KKT\_L2, for estimating the atmospheric extinction coefficient from Raman lidar data. The inverse problem under investigation is non--linear, and we have used the Karush--Kuhn--Tucker conditions to incorporate a non--negativity constraint and derive the optimization algorithms. The main novelty of the proposed methods is that they have been devised to model the Poisson statistics of the data, differently from all methods available in the literature on this topic. 

We have compared the proposed algorithms with the state of the art in the field, namely Richardson--Lucy and Tikhonov. The first one assumes that noise is Poisson on the log--transformed data, which is incorrect. The second one, which also works with the log--transformed data, assumes that noise is Gaussian and with constant variance; in order to account for an altitude--variable noise variance, we have also introduced a weighted Tikhonov approach, where the weights have been obtained from the empirical variance of the signals. 

The numerical results on synthetic and experimental data appear consistent and indicate that the proposed methods outperform the state of the art. In particular, KKT is better than its most direct competitor RL, inasmuch the estimated profiles are either less affected by noise or less over--smoothed. This improvement, clearly visible but not striking, was expected due to the incorrect assumption on the noise distribution in RL. 
The comparison of KKT\_L2 with Tikhonov is more striking: while the improvement over plain Tikhonov is obvious, because the noise variance changes strongly across the altitudes and plain Tikhonov fails to take this into account, the fairly substantial improvement over weighted Tikhonov comes a bit more unexpected. We speculate that the relatively unsatisfactory performances of weighted Tikhonov might be due either to a non--optimal choice of the weighting matrix, or quite possibly to the fact that the Gaussian distribution does not model correctly the log--transformed data, which are likely to be more skewed than a Gaussian can manage. However, more work is needed to better address this point.

Our results seem also to point out that the early stopping approach is, in this case, less effective than regularization through an explicit $\ell^2$ term. 
Indeed, while the profiles estimated by KKT are still affected by noise--induced spikes, or, alternatively, by some degree of oversmoothing, those obtained by KKT\_L2 present a homogeneous behaviour across altitudes. In practice, when using KKT it is almost necessary to tune the stopping iteration differently, depending on whether one is interested in the lower or higher part of the range. Instead, KKT\_L2 provides a reliable reconstruction at both low and high altitudes with the same value of the regularization parameter.

The results in this study further add to the field of inverse problems with Poisson noise, a lively field, attracting increasing interest recently \cite{hohage2016inverse}. We believe there are several interesting directions for future work. One is to consider an altitude--dependent regularization parameter by means of multi--parameter approaches, such as the one in \cite{bortolotti2016uniform}. Alternatively, it would be interesting to investigate the use of more sophisticated regularizers, such as Total Variation \cite{ng2010solving} or sparsity--inducing norms on a frame expansion, like has been done in \cite{melot2012some} for a problem with the same structure. Of course, due to the presence of non--smooth terms, in this case more sophisticated algorithms, such as proximal and primal--dual methods \cite{combettes2011proximal,chouzenoux2014variable, condat2013primal}, should be used.

\section*{Acknowledgments}\label{sec:sec6}
The authors acknowledge EARLINET for providing aerosol synthetic lidar profiles. The authors also acknowledge Advanced Lidar Applications srl for providing real lidar measurements.
AS has been partially supported by Gruppo Nazionale per il Calcolo Scientifico.
\section*{References}
\bibliography{biblio}
\end{document}